\documentclass{article}

\usepackage[english]{babel}

\usepackage{amsthm,amsmath}
\usepackage[letterpaper,top=2cm,bottom=2cm,left=3cm,right=3cm,marginparwidth=1.75cm]{geometry}
\theoremstyle{definition}
\newtheorem{definition}{Definition}
\theoremstyle{theorem}
\newtheorem{theorem}{Theorem}
\usepackage{amsmath}
\usepackage{graphicx}
\usepackage[colorlinks=true, allcolors=blue]{hyperref}
\usepackage{amsfonts}
\DeclareMathOperator{\rd}{d\!}

\title{On the uniqueness of Minimal Translation Surfacess in the Heisenberg Group}
\author{Christiam Figueroa}

\begin{document}
\maketitle

\begin{abstract}
In this paper, we establish a complete classification of minimal translation surfaces within the 3-dimensional Heisenberg group $\mathcal{H}_3$. While translation surfaces in $\mathcal{H}_3$ yield diverse families of solutions depending on the choice of generating curves, we prove that, up to isometries of $\mathcal{H}_3$, there exists a unique class of translation surfaces.
\end{abstract}

\section{Introduction}

Its known that the minimal translation surface in the Euclidena plane, besides the planes, are given by,
$$z(x,y)=\frac{1}{a}\ln|\frac{\cos ax}{\sin ax}|$$
A translation surface in any 3-dimensional Lie group $G$ equipped with a left-invariant Riemannian metric is a surface in $G$ parameterized as the product of two curves. In \cite{inoguchi2012minimal}, Inoguchi et al. gave a classification of minimal translation surfaces in the Heisenberg group obtained as the product of two curves situated in the coordinate planes of $\mathcal{H}_3$.In this paper, we prove that, up to isometries of $\mathcal{H}_3$, there exists a unique class of translation surfaces.

\section{The Gans Model }
This is a  model of the hyperbolic geometry, developed by David Gans, see \cite{gans1966new}.
Consider the Poincar\'{e} Disk
\[\mathbb{D}=\{(x,y): x^{2}+y^{2}<1\}\]
endowed with the metric
\[g(x,y)=\frac{4}{(1-x^{2}-y^{2})^{2}}(\rd x^{2}+\rd y^{2}). \]

We   define a  diffeomorphism between   the Poincar\'{e} disk and the plane $\mathcal{P}: z=1$. To do this we use the stereographic projection from the south pole $(0,0,-1)$ of the unit sphere, $S^2$, we define the following diffeomorphism  $\varphi$ between the upper hemisphere $S_{+}^{2}$ onto  the disk, $\mathbb{D}\subset \mathbb{R}^3$
\[\varphi(x,y,z)=(\frac{x}{z+1},\frac{y}{z+1},0).
\]
Similarly, considering  the stereographic projection from the origin $(0,0,0)$ of $S^2$, we define a diffeomorphism $\psi$ of $S_{+}^{2}$ onto  the plane  $\mathcal{P}: z=1$,
\begin{equation}\label{psi}
  \psi(x,y,z)=(\frac{x}{z},\frac{y}{z},1).
\end{equation}
Then, $F(x,y)=\psi\circ\varphi^{-1}$ is a diffeomorphism  from the disk $\mathbb{D}$ onto $\mathcal{P}$, where
\begin{equation}\label{diff}
  F(x,y,0)=(\frac{2x}{1-x^{2}-y^{2}},\frac{2y}{1-x^{2}-y^{2}},1)
\end{equation}
and the inverse is given by
\[F^{-1}(u,v,1)=(\frac{u}{1+\sqrt{1+u^{2}+v^{2}}},\frac{v}{1+\sqrt{1+u^{2}+v^{2}}},0)\]
Then the metric induced on  $\mathcal{P}$ by $F$ is given by
\[h(u,v)=\frac{(1+v^{2})\rd u^{2}-2uv\rd u\rd v+(1+u^{2})\rd v^{2}}{1+u^{2}+v^{2}}\]
The Riemannian space $(\mathcal{P},h)$ is the Gans model of  the hyperbolic geometry.

\subsection{Isometries}

Consider the Poincar\'{e} disk  as the subset $\mathbb{D}=\{z\in  \mathbb{C}: |z|<1\}$
of the complex plane and the Gans model $\mathcal{P}=\{w:w\in \mathbb{C}\}$.
We know that the set of orientation-preserving isometries of the Poincar\'{e} Disk have the form,
\[\rho(z)=e^{i\theta}\frac{z-a}{1-\overline{a}z}, \qquad a\in\mathbb{D}.\]
And all  isometries of $\mathbb{D}$ are composed of $\rho$ with complex conjugation, that is reflection at the real axis.
Therefore, the isometry group of the Gans model is
\[\mathrm{Iso}(\mathcal{P})=\{F\circ\rho\circ F^{-1}: \rho\in \mathrm{Iso}(\mathbb{D})\},\]
where $F$ is as in $(\ref{diff})$. I shall highlight two cases:

If $\rho(z)=e^{i\theta}z$ , then  $F\circ\rho\circ F^{-1}(w)=e^{i\theta}w$, that is, a rotation about the origin $(0,0)$ is an isometry of the hyperbolic space $\mathcal{P}$.

On the other hand, if $\rho(z)=\overline{z}$, then  $F\circ\rho\circ F^{-1}(w)=\overline{w}$  is the
 reflection across the $u$ axis. Since rotation about the origin is an isometry, a reflection across the line $au+bv=0$  is an isometry too.

\section{The Geometry of the Heisenberg Group }
The 3-dimensional Heisenberg group $\mathcal{H}_{3}$ is a two-step nilpotent Lie group. It has the following standard representation
 in $GL_{3}(\mathbb{R})$
 \[\left[
      \begin{array}{ccc}
        1 & r & t \\
        0 & 1 & s \\
        0 & 0 & 1 \\
      \end{array}
    \right]
\]
with $r,s,t\in \mathbb{R}$.

In order to describe a left-invariant metric on $\mathcal{H}_{3}$, we note that the Lie algebra $\mathfrak{h}_{3}$ of $\mathcal{H}_{3}$ is given by the
matrices
\[A=\left[
      \begin{array}{ccc}
        0 & x & z \\
        0 & 0 & y \\
        0 & 0 & 0 \\
      \end{array}
    \right]
\]
    with $x,y,z$ real. The exponential map $\exp:\mathfrak{h}_{3}\rightarrow \mathcal{H}_{3}$ is a global diffeomorphism, and is given by
    \[ \exp(A)=I+A+\frac{A^{2}}{2}=\left[
                                   \begin{array}{ccc}
                                     1 & x & z+\frac{xy}{2} \\
                                     0 & 1 & y \\
                                     0 & 0 & 1 \\
                                   \end{array}
                                 \right].
    \]
    Using the exponential map as a global parametrization, with the identification of the Lie algebra $\mathfrak{h}_{3}$ with $\mathbb{R}^{3}$ given by
\[(x,y,z)\longleftrightarrow \left[
      \begin{array}{ccc}
        0 & x & z \\
        0 & 0 & y \\
        0 & 0 & 0 \\
      \end{array}
    \right]
\]
    the group structure of $\mathcal{H}_{3}$ is given by
\begin{equation}\label{pr}
(a,b,c)\ast (x,y,z)=(a+x,b+y,c+z+\frac{ay-bx}{2}).
\end{equation}
    From now on, modulo the identification given by $exp$, we consider $\mathcal{H}_{3}$ as $\mathbb{R}^{3}$ with the product given in (\ref{pr}). The Lie   algebra bracket, in terms of the canonical basis $\{e_{1},e_{2},e_{3}\}$ of $\mathbb{R}^{3}$, is given by
    \vspace{0.5cm}
    \[[e_{1},e_{2}]= e_{3}, \qquad   [e_{i},e_{3}]=0,\]
    \vspace{0.5cm}
 with $i=1,2,3.$ Now, using $\{e_{1},e_{2},e_{3}\}$ as the orthonormal frame at the identity, we have the following left-invariant metric $\rd s^{2}$ in $\mathcal{H}_{3}$
\[\rd s^{2}=\rd x^{2}+\rd y^{2}+(\frac{1}{2}y \rd x-\frac{1}{2}x\rd y+\rd z)^{2}.\]
And the basis of the orthonormal left-invariant vector fields is given by
    \[E_{1}=\frac{\partial}{\partial x}-\frac{y}{2}\frac{\partial}{\partial z},\qquad  E_{2}=\frac{\partial}{\partial x}+\frac{x}{2}\frac{\partial}{\partial z},\qquad
    E_{3}=\frac{\partial}{\partial z}\cdot
    \]
Then the Riemann connection of $\rd s^{2}$, in terms of the basis $\{E_{i}\}$, is given by

\[
\begin{array}{ccccc}
  \nabla_{E_{1}}E_{2} & =& \frac{1}{2}E_{3}  &=&  -\nabla_{E_{2}}E_{1} \\\\
  \nabla_{E_{1}}E_{3} & =& -\frac{1}{2}E_{2} &=&  \nabla_{E_{3}}E_{1} \\\\
  \nabla_{E_{2}}E_{3} & =& \frac{1}{2}E_{1}  &=&  \nabla_{E_{3}}E_{2}
\end{array}
\]
and $\nabla_{E_{i}}E_{i}=0$ for $i=1,2,3$.

Using the fact that an isometry of $\mathcal{H}_{3}$ which fix the identity, is an automorphism of $\mathfrak{h}_{3}$, it is possible to show that every isometry of $\mathcal{H}_{3}$ is of the form $L\circ A$ where  $L$ is a left translation in $\mathcal{H}_{3}$ and $A$  is in one of the following forms
\[\begin{bmatrix}
    \cos\theta & -\sin\theta\;\;\, & 0 \\
    \sin\theta & \cos\theta &0 \\
    0 & 0 & 1 \\

  \end{bmatrix}
\qquad \textrm{or} \qquad  \begin{bmatrix}
     \cos\theta & \sin\theta & 0 \\
    \sin\theta & -\cos\theta \;\;\, &0 \\
    0 & 0 & -1\;\; \\
   \end{bmatrix}.
   \]

 That is, $A$ represent a rotation around the $z$-axis or a composition of the reflection across the plane $z=0$ and a reflection across a line $y=mx$ for some $m\in \mathbb{R}$ .

 \section{Surfaces in the Heisenberg Group}

Let $S$ be a graph of a smooth function $f:\Omega \rightarrow \mathbb{R}$ where $\Omega $ is an open set of $\mathbb{R}^{2}$. We consider the
following parametrization of $S$
\begin{equation}\label{paramet}
 X\left( x,y\right) =( x,y,f( x,y)),\qquad (x,y)\in \Omega.
\end{equation}
A basis of the tangent space $T_{p}S$ associated to this
parametrization is given by
\begin{equation}
\begin{array}{ccccc}
X_{x} & = & \left( 1,0,f_{x}\right) & = & E_{1}+\left( f_{x}+\frac{y}{2}%
\right) E_{3} \\\\
X_{y} & = & \left( 0,1,f_{y}\right) & = & E_{2}+\left( f_{y}-\frac{x}{2}%
\right) E_{3}%
\end{array}
\label{basis}
\end{equation}%
\noindent
where $f_x$ and $f_y$ denote the partial derivatives of $f$, with respecto $x$ and $y$ respectvely. And the  unit normal vector of $S$  is given by
\begin{equation}
\eta \big( x,y\big) =-\bigg(\frac{f_{x}+\displaystyle\frac{y}{2}}{w}\bigg)
E_{1}-\bigg( \frac{f_{y}-\displaystyle\frac{x}{2}}{w}\bigg) E_{2}+\frac{1}{w}E_{3}
\label{normal}
\end{equation}%
where
\begin{equation}
w=\sqrt{1+\left( f_{x}+ \frac{y}{2}\right) ^{2}+\left( f_{y}-\frac{x}{2}\right) ^{2}}.
\label{w}
\end{equation}
Then the
coefficients of the first fundamental form of $S$  are given by%
\begin{equation}
\begin{array}{ccccl}
E & = & \langle X_{x}, X_{x}\rangle & = & 1+\left( f_{x}+ \displaystyle\frac{y}{2}\right) ^{2} \\\\
F & = & \langle X_{y},X_{x}\rangle & = & \left( f_{x}+\displaystyle\frac{y}{2}\right) \left( f_{y}-%
\displaystyle\frac{x}{2}\right) \\\\\
G & = & \langle X_{y},X_{y}\rangle & = & 1+\left( f_{y}- \displaystyle\frac{x}{2}\right) ^{2}.%
\end{array}
\label{1ffund}
\end{equation}%
If $\nabla $ is the Riemannian connection of $\left( \mathcal{H}%
_{3},\rd s^{2}\right) $, by the Weingarten  formula for hypersurfaces, we have that%
\[
A_{\eta }v=-\nabla _{v}\eta ,\qquad  v\in T_{p}S
\]%
and the coefficients of the second fundamental form are given by
\begin{equation}
\begin{array}{ccccl}
L & = & -\langle\nabla _{X_{x}}\eta ,X_{x}\rangle & = & \displaystyle\frac{f_{xx}+( f_{y}-\frac{x%
}{2})( f_{x}+\frac{y}{2})}{w} \\\\
M & = & -\langle\nabla _{X_{x}}\eta ,X_{y}\rangle & = &\displaystyle \frac{f_{xy}+\frac{1}{2}\left(
f_{y}-\frac{x}{2}\right) ^{2}-\frac{1}{2}\left( f_{x}+\frac{y}{2}\right) ^{2}%
}{w} \\\\
N & = & -\langle\nabla _{X_{y}}\eta ,X_{y}\rangle & = & \displaystyle\frac{f_{yy}-\left( f_{y}-\frac{x%
}{2}\right) \left( f_{x}+\frac{y}{2}\right) }{w}\cdot%
\end{array}
\label{2ffund}
\end{equation}

Recall that the mean curvature of any surface of $\mathcal{H}_3$ can be expressed in terms of its first and second fundamental forms, given. a parametrization,
\[H=\frac{1}{2}\big(\frac{EN+GL-2FM}{EG-F^2}\big)\]
Given a parametrization when the surface is the graph of a smooth function $f$, we replace the coefficients  into the mean curvature formula,
\[\frac{(1+q^2)f_{xx}-2pqf_{xy}+(1+p^2)f_{yy}}{(1+p^2+q^2)^{3/2}}=2H,\]
where $p=f_x+y/2$ and $q=f_y-x/2$. In particular, when $H=0$ the equation of the minimal graph is given by
\begin{equation}\label{meq}
  (1+q^2)f_{xx}-2pqf_{xy}+(1+p^2)f_{yy}=0
\end{equation}

\section{The Gauss Map}
Recall that the Gauss map is a function  from an oriented surface, $S\subset \mathbb{E}^{3}$, to the unit sphere in the Euclidean space . It associates to every point on the surface its oriented unit normal vector. Considering the Euclidean space as a commutative Lie group, the Gauss map is just the translation of the unit normal vector at any point of the surface to the origin, the identity element of $\mathbb{R}^{3}$. Reasoning in this way we define a Gauss map in the following form

\begin{definition}
Let $S\subset G$ be an orientable hypersurface of a n-dimensional Lie group $%
G,$ provided with a left invariant metric. The map%
\[
\gamma :S\rightarrow S^{n-1}=\left\{ v\in \tilde{g}: \left\vert v\right\vert
=1\right\}
\]
where $\gamma \left( p\right) =\rd L_{p}^{-1}\circ \eta \left( p\right) $, $\tilde{g}$ the Lie algebra of $G$
and $\eta $ the unitary normal vector field of $S,$ is called the Gauss map of $S.$
\end{definition}

\noindent We observe that
\[\rd\gamma\left( T_{p}S\right) \subseteq T_{\gamma \left( p\right) }S^{n-1}
 =  \left\{ \gamma \left( p\right) \right\} ^{\perp }  =
\rd L_{p}^{-1}\left( T_{p}S\right),\]

\noindent therefore $\rd L_{p}\circ \rd\gamma \left( T_{p}S\right)
\subseteq T_{p}S$ .

Now we  obtain a local expression of the Gauss map $\gamma $. In fact, we consider the following sequence of maps
\[\phi:\Omega \overset{X}\longrightarrow X(\Omega)\subset \mathcal{H}_{3}\overset{\gamma}\longrightarrow S^{2}\overset{\psi}\longrightarrow \mathcal{P}\]
where, $X$ is a parametrization of $S$, $\mathcal{P}$ is the plane $z=1$ and $\psi$ is given by $(\ref{psi})$.

In the case where $S$ is the graph of a smooth function $f\left(x,y\right)$ with $(x,y)$ over a domain $\Omega \subset \mathbb{R}^{2} $, the Gauss map is given by

\begin{equation}\label{gmap}
  \phi(x,y)=\left(-(f_{x}+\frac{y}{2}),-(f_{y}-\frac{x}{2})\right)
\end{equation}
and the Jacobian matrix  of $\phi$ is
\begin{equation}\label{jacob}
\rd \phi_{(x,y)}=\left(
                   \begin{array}{cc}
                     -f_{xx} & -f_{xy}-1/2 \\
                     -f_{xy}+1/2 & -f_{yy} \\
                   \end{array}
                 \right).
\end{equation}

\noindent The determinant of $\rd\phi$ at a point $(x,y)$ is given by
\begin{equation}
\det  \rd \phi_{(x,y)} =f_{xx}f_{yy}-f_{xy}^{2}+\frac{1}{4} \label{rank}
\end{equation}%
 which we define as the determinant of the Gauss map at  $(x,y).$ In particular, when $\phi$ is non-constant and its determinant vanishes identically, we say that the Gauss map has rank 1.

  If $\Omega=\mathbb{R}^{2}$, the greatest lower bound of the absolute value
of $\det\rd \phi_{(x,y)}$ is zero. This was proved by A. Borisenko and E.  Petrov in \cite{borisenko2011surfaces}.

We know that in the Euclidean case the differential of the Gauss
map is just the second fundamental form for surfaces in
$\mathbb{R}^{3},$ and this fact can be generalized for hypersurfaces
in any Lie group. The following theorem, see \cite{ripoll1991hypersurfaces},
states a relationship between the Gauss map and the extrinsic
geometry of $S.$

\begin{theorem}
\label{gauss}Let $S$ be an orientable hypersurface of a Lie group. Then
\[
\rd L_{p}\circ \rd\gamma _{p}\left( v\right) =-\left( A_{\eta }\left( v\right)
+\alpha _{\bar{\eta}}\left( v\right) \right) ,\qquad v\in T_{p}S
\]
where $A_{\eta }$ is the Weingarten operator, $\alpha _{\bar{\eta}}\left(
v\right) =\nabla _{v}\bar{\eta}$ and $\bar{\eta}$ is the left invariant
vector field such that $\eta \left( p\right) =\bar{\eta}\left( p\right).$
\end{theorem}

In \cite{figueroa2012gauss}, several results were obtained regarding minimal surfaces in the Heisenberg group and their relation to the Gauss map.
\begin{theorem}
\label{vertical}There is no graph of a smooth function over $XY$ with constant Gauss map. More over, vertical planes are the only minimal surfaces whose Gauss map is constant.
\end{theorem}

 \begin{theorem}{\label{rank1}}
   The minimal graphs in $\mathcal{H}_3$ with Gauss map of rank one, are

 \begin{equation}\label{gramin}
   f(x,y)=\frac{xy}{2}-\frac{k}{2}[y\sqrt{1+y^2}+\ln(y+\sqrt{1+y^2})]
 \end{equation}
   with $k\in \mathbb{R}$
 \end{theorem}

\section{Minimal translation surfaces}
In \cite{inoguchi2012minimal}, Inoguchi define and study translation surface in $\mathcal{H}3$. A translation surface $M(\gamma_1,\gamma_2)$ is a surface parameterized by
$$r(x,y)=\gamma_1(x)*\gamma_2(y)$$
where $\gamma_1$ and $\gamma_2(y)$ are curves situated in the coordinates planes of $\mathbb{R}^3$. Since the operation $*$ is not commutative, they distinguish six types of minimal translation surfaces in $\mathcal{H}_3$.

\subsection{Surfaces type 1 and 4}
In this case the curves are $\gamma_1(x)=(x,0,u(x))$ and $\gamma_2(y)=(0,y,v(y))$. The translation surfaces are parameterized by

\begin{eqnarray}
  r_1(x,y) &=& \gamma_1(x)*\gamma_2(y)=(x,y,u(x)+v(y)+\frac{xy}{2} \\
  r_4(x,y) &=&  \gamma_2(x)*\gamma_1(y)=(x,y,u(x)+v(y)-\frac{xy}{2}
\end{eqnarray}
We notice that in both cases, these surfaces are graph of a function and both are isometric through a $90°$ rotation. From (\ref{gmap}),  the  Gauss map of $r_1$  is

$$\phi(x,y)=\left(-(f_{x}+\frac{y}{2}),-(f_{y}-\frac{x}{2})\right)=(u'(x)+y,v'(y))$$
and its matrix jacobian is,
$$\rd \phi_{(x,y)}=\left(
                           \begin{array}{cc}
                             u''(x) & -1 \\
                             0 &  v''(y)
                           \end{array}
                         \right)
                 $$
Inoguchi proved, for the minimal case,  it must be true that $u''(x)=0$ or $v''(y)=0$, that is, the rank of its Gauss map must be 1. According to theorem ($\ref{rank1}$) and up to isometries, these two types of surface are graph of the following function
$$f(x,y)=\frac{xy}{2}-\frac{k}{2}[y\sqrt{1+y^2}+\ln(y+\sqrt{1+y^2})],
$$
where $k\in \mathbb{R}$.

\subsection{Surfaces type 2 and 5}
In this case the curves are $\gamma_1(x)=(x,0,u(x))$ and $\gamma_2(y)=(v(y),y,0)$. The translation surfaces are parameterized by

\begin{eqnarray}
  r_2(x,y) &=& \gamma_1(x)*\gamma_2(y)=(x+v(y),y,u(x)+\frac{xy}{2}) \\
  r_5(x,y) &=&  \gamma_2(x)*\gamma_1(y)=(x+v(y),y,u(x)-\frac{xy}{2})
\end{eqnarray}
For $r_2$ the normal direcction to this surface, in terms of the left invariant fields, is given by
$$N_2=-\Big(y+u'(x)\Big)E_1+\Big(u'(x)v'(y)+\frac{1}{2}(yv'(y)+v(y))\Big)E_2+E_3$$
So the Gauss map representation given by (\ref{psi}) is
$$\phi_2(x,y)=\Big(y+u'(x),u'(x)v'(y)+\frac{1}{2}(yv'(y)+v(y))\Big)$$
The Jacobian matrix of its Gauss map is
$$J\phi_2(x,y)=\left(
                 \begin{array}{cc}
                   -u''(x) & -1 \\
                   u''(x)v'(y) & u'(x)v''(y)+v'(y)+yv''(y)/2 \\
                 \end{array}
               \right)
$$
and its determinant
$$\det J\phi_2(x,y)=-u''(x)v''(y)(\frac{y}{2}+u'(x)).$$
As in the previous case  according to  Inoguchi, $u''(x)=0$ or $v''(y)=0$, hence this surface has rank 1. More over, it´s easy to verify that $r_2$ is injective and that every straight line parallel to the $Z$-axis intersects this surface at most at one point, therefore, this surface can be parameterized as the graph of a function over the $XY$-plane.

For surface type 5, we proceed in the same way to reach the same conclusion.

\subsection{Surfaces type 3 and 6}
In this case the curves are $\gamma_1(x)=(0,x,u(x))$ and $\gamma_2(y)=(y,v(y),0)$. The translation surfaces are parameterized by

\begin{eqnarray}
  r_3(x,y) &=& \gamma_1(x)*\gamma_2(y)=(y,v(y)+x,u(x)-\frac{xy}{2} \\
  r_6(x,y) &=&  \gamma_2(x)*\gamma_1(y)=(y,v(y)+x,u(x)+\frac{xy}{2}
\end{eqnarray}
For $r_3$ the normal direcction to this surface, in terms of the left invariant fields, is given by
$$N_3=-\Big(\frac{1}{2}(yv'(y)+v(y))-u'(x)v'(y)\Big)E_1+\Big(y-u'(x)\Big)E_2-1E_3$$
So the Gauss map representation given by (\ref{psi}) is
$$\phi_3(x,y)=\Big((\frac{1}{2}(yv'(y)+v(y))-u'(x)v'(y)),u'(x) -y\Big)$$
The Jacobian matrix of its Gauss map is
$$J\phi_3(x,y)=\left(
                 \begin{array}{cc}
                   -v'(y)u''(x) & v'(y)+\frac{1}{2}yv''(y)-v''(y)u'(x) \\
                   u''(x) & -1 \\
                 \end{array}
               \right)
$$
and its determinant
$$\det J\phi_3(x,y)=u''(x)v''(y)(u'(x)-\frac{y}{2})$$
As is the previous case and according to  Inoguchi, $u''(x)=0$ or $v''(y)=0$, so this surface has rank 1 and it´s also easy to verify that $r_3$ is injective and every straight line parallel to the $Z$-axis intersects this surface at most at one point, so it can also be parameterized as the graph of a function over the $XY$-plane.
For surface type 6, we proceed in the same way to reach the same conclusion.

All the  aforementioned cases can be summarized in the following theorem

\begin{theorem}
  If $S$ is a minimal translation surface in the Heisenberg group, up to isometries, this surface is given by the graph of the following function:
  $$f(x,y)=\frac{xy}{2}-\frac{k}{2}[y\sqrt{1+y^2}+\ln(y+\sqrt{1+y^2})]$$
  where $k\in \mathbb{R}$
\end{theorem}

\end{document}